\documentclass[12pt]{article}
\usepackage{amsmath,amsfonts,amssymb,amsthm}

\newtheorem{proposition}{Proposition}[section]
\newtheorem{theorem}[proposition]{Theorem}

\newtheorem{lemma}[proposition]{Lemma}
\newtheorem{rem}[proposition]{Remark}
\newtheorem{cor}[proposition]{Corollary}
\newtheorem{Def}[proposition]{Definition}

\newcommand{\dd}{{\mathrm d}}
\newcommand{\e}{{\mathrm e}}
\newcommand{\Oh}{{\mathrm O}}
\newcommand{\oh}{{\mathrm o}}

\newcommand{\Prob}{{\mathbb P}}
\newcommand{\Probx}{{\mathbb P}^{(x)}}
\newcommand{\Exp}{{\mathbb E}}

\newcommand{\Gb}{\overline G}
\newcommand{\Fb}{\overline F}

\newcommand{\nc}{\newcommand}

\nc{\I}{{\mathbf 1}}
\nc{\R}{{\mathbb R}}
\nc{\N}{{\mathbb N}}
\nc{\Z}{{\mathbb Z}}

\nc{\BP}{\mathbb{P}}
\nc{\BE}{\mathbb{E}}
\nc{\BQ}{\mathbb{Q}}
\nc{\BX}{\mathbb{X}}
\numberwithin{equation}{section}
\newcommand{\halmos}{{\mbox{\, \vspace{3mm}}} \hfill
\mbox{$\Box$}}

\begin{document}
\renewcommand{\thefootnote}{\fnsymbol{footnote}}

\author{S\o ren Asmussen,\thanks{Postal adress: Department of Mathematics, Aarhus University, Ny Munkegade, DK-8000 Aarhus C, Denmark. Email address: asmus@imf.au.dk} {\em \normalsize Aarhus University} \and Sergey
Foss,\thanks{Postal adress: School of Mathematical and Computer Sciences, Heriot-Watt University, EH14 4AS, Edinburgh, United Kingdom.
Email address: s.foss@hw.ac.uk. This author thanks MSRI, Berkeley, for their hospitality.} {\em \normalsize Heriot-Watt University and Institute of Mathematics, Novosibirsk}}
\title{ON EXCEEDANCE TIMES\\ FOR SOME PROCESSES\\ WITH DEPENDENT INCREMENTS}
\date{\today}
\maketitle

\newcommand{\rtreg}{\tau}
\newcommand{\rtrw}{\tau^{\rm rw}}
\newcommand{\rpreg}{\psi}

\begin{abstract}
Let $\{ Z_n\}_{n\ge 0}$ be a random walk with a negative drift and
i.i.d. increments with heavy-tailed distribution and let
$M=\sup_{n\ge 0}Z_n$ be its supremum.  Asmussen \& Kl{\"u}ppelberg
(1996) considered the behavior of the random walk given that $M>x$,
for $x$ large, and obtained a limit theorem, as $x\to\infty$, for
the distribution of the quadruple that includes the time $\rtreg=\rtreg(x)$ to exceed
level $x$, position $Z_{\rtreg}$ at this time, position $Z_{\rtreg-1}$ at the prior time, and
the trajectory up to it (similar results were obtained for the
Cram\'er-Lundberg insurance risk process). We obtain here several extensions of this
result to various regenerative-type models and, in particular,
to the case of a random walk with dependent increments. Particular attention is given
to describing the limiting conditional behavior of $\tau$.
The class of models include Markov-modulated 
models as particular cases. 
We also study fluid models, the Bj{\"o}rk-Grandell risk process,  give examples
where the order of $\tau$ is genuinely different from the random walk case,
and discuss which growth rates are possible.
Our proofs are purely probabilistic and are
based on results and ideas from Asmussen, Schmidli \& Schmidt (1999), Foss
\& Zachary (2002), and Foss, Konstantopoulos \& Zachary (2007).\\[2mm]
\end{abstract}

\newpage

\section{Introduction}\label{S:Intr}

Let $Z=\{Z(t)\}_{t\ge 0}$ be a stochastic process with increments having a regenerative structure (\cite{APQ}):
there exist random times $T_0=0,T_1,T_2,\ldots$  
splitting $Z$ up into i.i.d.\ cycles
$$\bigl\{Z(t)-Z(0)\bigr\}_{0\le t<R_1}\,=\,\bigl\{Z(t+T_0)-Z(T_0)\bigr\}_{0\le t<
R_1},\ \ \bigl\{Z(t+T_k)-Z(T_k)\bigr\}_{0\le t<R_{k+1}},\
\ldots
$$
with lengths $R_0=T_0=0,R_1=T_1-T_0,R_2=T_2-T_1,\ldots$
(traditionally as in \cite{APQ}, one allows the first cycle to have a different distribution; we won't do this
since our results
are easily adapted to this setting). We will also assume $Z(0)=0$.  
A main example we have in mind is the claims surplus process of an
insurance company (accumulated claims minus premiums, cf.\ \cite{RP2}). In this setting,
$\rtreg\, =\, \rtreg(x)\, =\, \inf\{t:\, Z(t)>x\}$ is the ruin time with initial surplus $x$, $M=\sup_{t\ge 0} Z(t)$
is the maximal claims surplus, and
$$\Prob(\rtreg<\infty)\ =\ \Prob(M>x)$$
is the ruin probability,
but $\rtreg$ and $M$ are also of interest
in many other contexts. For example, $M$  could be the stationary waiting time in a
single-server queue with i.i.d.\ service times whose input process
is modulated by a Markov chain (say, this is an output process
from another stationary single-server queue, see e.g.
\cite{FBSF}). 

Inder suitable 
assumptions,  
the supremum $$
\sup_{0\le t\le R_{k+1}}(Z(t+T_k)-Z(T_k))
$$ over a typical regenerative cycle of the process increments
has a heavy-tailed distribution, say $F$, on $[0,\infty)$ with mean $m_F<\infty$ whose
{\it integrated tail} distribution
$${\overline{F^I}}(x)\,=\,\min \left(
1, \int_{x}^{\infty}\overline{F}(y)\,\dd y\right)$$
is subexponential. Then 
it has been proved in a variety
of settings  that
\begin{equation}\label{S15.6b0}
\Prob(\rtreg<\infty)\ =\ 
\Prob(M>x)\ \sim\ b{\overline{F^I}}(x)\,,\ \ x\to\infty,
\end{equation}
where
$b>0$ is a constant, thereby extending a classical result for
random walks and the Cram\'er-Lundberg process due to (in
alphabetical order) Borovkov, Cohen, Embrechts, Pakes, Veraverbeke, von
Bahr. In particular, Asmussen, Schmidli \& Schmidt~\cite{ASS}
proved the following (for background on subexponential
distributions, see, e.g.,  \cite{EmKluMik}, \cite[X.1]{RP2}, or
\cite{SFDKSZ}):

\begin{theorem}\label{Th_ASS} In the regenerative setting, let
$$\xi_k\ =\ Z(T_{k+1})-Z(T_k)\,,\ \ \
\xi^*_k\ =\ \sup_{T_k\le t<T_{k+1}}Z(t)-Z(T_k)\,.$$
Assume that
\begin{equation}\label{S15.6b}\Prob(\xi_1>x)\,\sim\,\Prob(\xi_1^*>x)\,\sim\,\Fb(x)\,,\ \ x\to\infty,
\end{equation}
for some distribution $F$ such that 
${\overline{F^I}}$ is a subexponential tail and that
$-a=\Exp\xi_1<0$.
Then $$\Prob(M>x)\,\sim\,
\frac{1}{a}{\overline{F^I}}(x),\quad x\to\infty.$$
\end{theorem}
As demonstrated by the examples in \cite{ASS} (and later papers, of which Asmussen \& Biard~%
\cite{SARB} is a recent instance),
this result covers a large number of examples.
Foss \& Zachary~\cite{SFSZ} gave a similar result in the case of 
a modulated random walk.

The purpose of the present paper is to supplement Theorem~\ref{Th_ASS} and the 
corresponding result from \cite{SFSZ} with
a description of the asymptotic behavior of $\rtreg$ given $M>x$, but in a more
general setting that covers both scenarios (of regenerative structure and of 
modulation). Results of this type
were given for the first time in Asmussen \& Kl{\"u}ppelberg~\cite{AK}, assuming that
$Z$ is either the classical Cram\'er-Lundberg risk process, a L\'evy process, or a discrete time
random walk $Z_n=\xi_1+\cdots+\xi_n$ with the $\xi_k$ i.i.d.\ and having common distribution
$F$ and mean $-a<0$. Note that there is a  discrete time
random walk imbedded in the regenerative setting: consider the 
process $Z$ at times $T_n$.

In the random walk setting, the basic assumption of \cite{AK}
is that there exists a function $e(x)\uparrow\infty$
such that, for any $t>0$, 
\begin{equation}\label{conv}
\lim_{x\to\infty}
\frac{\overline{F^I}\bigl(x+te(x)\bigr)}{\overline{F^I}(x)}
\ =\  \Gb(t)
\end{equation}
for some distribution $G$. We assume in addition that the function 
$e(x)$ is what could be called {\it weakly self-neglecting}, i.e.
\begin{equation}\label{weakSN}
\limsup_{x\to\infty}
\frac{e\bigl(x+e(x)\bigr)}{e(x)}<\infty.
\end{equation}
Both assumptions \eqref{conv} and \eqref{weakSN} 
hold in the standard examples of subexponential distributions,
see e.g. \cite{GoRe} and \cite{BadH} for further details. In the regularly varying case $\Fb(x)=L(x)/x^\alpha$, a natural
scaling is $e(x)=x$; then \eqref{weakSN} is automatic
and $G$ is Pareto with $\Gb(t)=(1+t)^{-\alpha}$. 
For other subexponential distributions such as the lognormal and the heavy-tailed Weibull,
one may take $e(x)=\overline{F^I}(x)/\overline{F}(x)$ and then $G$ is standard exponential.
Let $W$ be a r.v.\ with distribution $G$.
Then, with $\rtrw(x)\,=\,\inf\{n:\, Z_n>x\}$,
it is shown in \cite{AK} (for later contributions in the same direction,
see \cite{klukypmal04}, 
\cite{SADK}) that:
\begin{theorem}\label{cor1}
Given in the random walk setting that $F^I$ is a subexponential distribution and that \eqref{conv} holds, 
as $x\to\infty$, the conditional distribution of
$\rtrw(x)/e(x)$ given $M>x$ converges to the distribution of $W/a$.
\end{theorem}
Our first main result
is the following extension. For a stochastic process with regenerative structure introduced
earlier, for cycle $i$, let 
$$
t_i \,= \,t_i(x) \,= \,\inf \{ t\le R_i :\, \ Z(t+T_{i-1})- Z(T_{i-1}) > x\} 
$$
if $\xi^{*}_i>x$, and $t_i=R_i$, otherwise.
\begin{theorem}\label{AFMain}
In the regenerative setting, 
assume in addition to the conditions
of Theorem~{\rm \ref{Th_ASS}} and to conditions \eqref{conv}-\eqref{weakSN} that for any $y>0$
\begin{equation}\label{S2506c}
\Prob\bigl(t_1>ye(x)\,\big|\,\xi_1>x\bigr)\ =\ \oh(1)\,,\ \ \ x\to\infty.
\end{equation}
Then the conditional distribution of
$\rtreg/e(x)$ given $M>x$ converges to the distribution of $\mu W/a$
where $\mu=\Exp R$.
\end{theorem}

The intuition behind Theorem~\ref{AFMain} is the following. In \cite{AK},
a number of supplementary results are given supporting the folklore that
exceedance of level $x$ occurs as result of one big $\xi_k$ and that 
all the other $\xi_k$ are `typical'. In the regenerative setting, it is shown in \cite{ASS} that the events
$\rtreg<\infty$ and $\rtrw<\infty$ (where the random walk is the process observed at times $T_n$) essentially are equivalent, and that exceedance asymptotically
occurs in cycle $\rtrw(x)$. Thus one expects  by the LLN, by the `typical' behavior before $\rtrw$
and by \eqref{S2506c} (which ensures that the length of the cycle in which ruin occurs
can be neglected), that conditionally on $M>x$, $\rtreg/\rtrw\to\mu$.
Given this, Theorem~\ref{cor1} then gives the desired conclusion.

The technical problem is to make this intuition precise in this
 and in more general settings. A difficulty is that conditioning on $\rtreg$ introduces
 some (presumably) small dependence between cycles $1,\ldots,\rtrw-1$ as well as some bias
 in their distribution (expected to be small as well); this was realized in \cite{SARB}, with the
 consequence that some results there are heuristic. To overcome this difficulty, we present
 an approach to results of type Theorem~\ref{cor1} 
 which is novel and combines the ideas from \cite{AK} and a sample-path analysis developed in
\cite{FBSF, SFTKSZ, SFSZ}.  The new approach is developed in Section~\ref{S:RW} in the setting of random walks
 modulated by a regenerative process $Y$. For such a process, the  asymptotics for $\Prob(M>x)$
 is given in \cite{SFTKSZ} (note that the setting allows $Y$ to be a Markov process with a
 general state space, whereas \cite{ASS} only can deal with the finite case). We supplement
 here with our second main result, Theorem 3.5, 
 giving the 
conditional behaviour of $\rtreg$. Compared to Theorem~\ref{AFMain}, it has the advantage 
that no conditions like  \eqref{S15.6b} or \eqref{S2506c}
 have to be verified, but it is also somewhat less general. 

It is easy to construct examples where \eqref{S2506c} fails as well as the conclusion of
Theorem~\ref{AFMain}, see Section~\ref{S:NSGrowth}. The order of $\rtreg$ may remain $e(x)$ (then with a larger multiplier than $\mu W/a$) or be effectively larger. It is tempting to conjecture that any rate
$\varphi (x)$ with $\varphi (x)/e(x)\to \infty$ may be attained. However, we shall show that $1/\Fb(x)$ is a critical upper bound.

\section{Preliminaries}\label{S:Prel}
We need some notation.

\begin{Def} Let $F$ be a distribution function and $\overline{F}(x)=1-F(x)$
its tail. Let  $h(x)$ be a positive non-decreasing function.
We say that $F$ is $h$-{\em insensitive} if
$$
\overline{F}\bigl(x+h(x)\bigr)\sim \overline{F}(x), \quad x\to\infty.
$$
\end{Def}
\noindent If \eqref{conv} holds for $F$, one can take $h$  as any function with
$h(x)=\oh\bigl(e(x)\bigr)$. One can find
more about the $h$-insensitivity property in \cite{SFDKSZ}, Chapter 2.
The term $h$-{\em flat} is also used by some authors, see e.g. \cite{balklures93}.

\begin{rem}\label{rem11} \rm Any subexponential distribution 
$F$ is long-tailed, i.e.\ 
$\overline{F}(x+C)\sim \overline{F}(x)$, for any constant $C$.
Therefore, by the diagonal argument, 
one  can choose a positive
function $h\uparrow\infty$ such that $F$ is also $h$-insensitive (clearly, the
choice of $h$ depends on $F$). If $F$ is $h$-insensitive and if
$0\le g \le h$, then $F$ is also $g$-insensitive.
\end{rem}

\begin{Def} We say that two families of events $A_{x}$ and $B_{x}$ of
positive probabilities, indexed by $x>0$, are {\it equivalent} and
write $A_{x}\sim B_{x}$, if $\Prob\bigl(A_x\Delta
B_x)\,=\,\oh\bigl(\Prob(A_x)\bigr)$, $x\to\infty$, where $A\Delta B=A\!\setminus\! B\,\cup\,
B\!\setminus\! A$ is the symmetric difference. 
\end{Def}
\noindent Note that if $A_{x}\sim B_{x}$, then also
$\Prob(A_x)\sim \Prob(B_x)$.

\section{Modulated random walk}\label{S:RW}

Consider a discrete-time regenerative process $Y= \{Y_n, n\ge 1\}$
such that, for each $n$, $Y_n$ takes values in some measurable
space $({\cal Y}, {\cal B_Y})$. We say that a random walk $\{Z_n,
n\ge 0\}$ defined by $Z_0=0$ and $Z_n = \xi_1+\cdots +\xi_n$ for
$n\ge 1$, is {\it modulated} by the process $Y$ if
\begin{description}
\item[{\rm (i)}]
conditionally on $Y$, the random variables $\xi_n, n\ge 1$, are
independent;
\item[{\rm (ii)}]
for some family $\{F_y, y\in {\cal Y}\}$ of distribution functions
such that, for each $x$, $F_y(x)$ is a measurable function of $y$,
we have, for $n=1,2,\ldots$,
\begin{equation}\label{cond11}
\Prob(\xi_n\le x \ | \ Y) = \Prob(\xi_n \le x \ | \
Y_n) = F_{Y_n}(x) \quad \mbox{a.s.}
\end{equation}
\end{description} 
Let $M^{rw}=\sup_{n\ge 0}Z_n$.
Under the conditions we give below,
$Z_n\to -\infty$ a.s. as $n\to\infty$, and so the random variable $M^{rw}$
is finite a.s.

The regenerative epochs of the modulating process $Y$ are denoted
by $0= T_0 < T_1 < \ldots $, 
with $R_k=T_k-T_{k-1}$. By definition, the cycles $\bigl(R_k,
(Y_n,0<n\le T_k-T_{k-1})\bigr), k\ge 1$, are i.i.d. We assume that
\begin{equation}\label{reg1}
\mu \,=\,\Exp R_1 <\infty.
\end{equation}
Let 
$$
\pi(B) \ = \ \frac{\Exp \sum_1^{R_1} {\mathbf 1} (Y_n\in
B)}{ \mu}, \quad B\in {\cal B_Y}
$$
be the stationary probability measure. We assume that each
distribution $F_y, y\in {\cal Y}$ has a finite mean
\begin{equation}\label{mean1}
a_y \ =\  \Exp[\xi_n\,|\, Y_n=y]\  =\ \int_{-\infty}^{\infty} x\, F_y(\dd x) \in (-\infty , \infty),
\end{equation}
and that 
\begin{equation}\label{ass1} \mbox{the family of distributions\ \ }
\{F_y , y\in {\cal Y}\} \   \mbox{ is uniformly
integrable}.
\end{equation}
In addition, we assume that this family of distributions satisfies
the following additional assumptions with respect to some 
reference distribution $F$ with finite mean and some measurable
function $c: {\cal Y} \rightarrow [0,1]$:
\begin{description}
\item[{\rm (C1)}]
$ \overline{F_y}(x) \le \overline{F}(x),$ for all $ x \in \R, y\in {\cal Y}$,
\item[{\rm (C2)}]
$\overline{F_y}(x) \sim c(y) \overline{F}(x)$ as  $x\to\infty$,
for all $y\in{\cal Y}$,
\item[{\rm (C3)}]
$\kappa \, =\, \sup_{y\in {\cal Y}}a_y$ is finite and $a\,=\, -
\int_{\cal Y} a_y \pi (\dd y)$ is finite and strictly positive,
\item[{\rm (C4)}]
for some nonnegative $b>\kappa$,
\begin{eqnarray*}
\Prob(bR_1>n) = \oh\bigl(\overline{F}(n)\bigr), \quad n\to\infty.
\end{eqnarray*}
\end{description}
Note that condition (C4) is redundant if $\kappa <0$ --- then one
can take $b=0$.

The following result is known (see Theorem 2.2 from \cite{SFTKSZ}
for a slightly more general version and also for discussion on
importance of conditions; see also Proposition 3.2 of \cite{klukypmal04}).
\begin{theorem}\label{Tnew2}
Suppose that conditions \eqref{cond11}--\eqref{ass1} and {\rm (C1)--(C4)}
hold and that the distribution $F^I$ is subexponential. Then $Z_n/n\to
-a $ a.s.\ as $n\to\infty$; in particular, $M^{rw}$ is an a.s.\ finite random
variable. Furthermore,
\begin{equation}\label{FKZ1}
\lim_{x\to\infty} \frac{\Prob(M^{rw}>x)}{{\overline{F^I}}(x)} =
\frac{C}{a}
\end{equation}
where $\displaystyle C\,=\,\int_{\cal Y} c(y) \pi (\dd y)\,\in\, [0,1]$ .
\end{theorem}
The main idea in the proof of Theorem \ref{Tnew2} is 
that 
the supremum of the modulated random walk, $M^{rw}$, may be closely approximated by a sum of two
{\it independent} random variables where one of them has a light-tailed distribution and the
other is the supremum of an ordinary random walk with i.i.d.
heavy-tailed increments with integrated tail distribution
proportional to $\overline{F}^I$.  

We also note that, by the strong law of large numbers (SLLN) and by the
diagonal argument, one can choose a sequence $\varepsilon_n\downarrow 0$
  such that
  \begin{equation}\label{SLL1}
  {\mathbb P} \bigl(|Z_{m}+ma|\le m\varepsilon_{m} \ \ \forall m\ge n\bigr)\to 1,
  \quad n\to\infty.
  \end{equation}
  Then
  \begin{equation}\label{SLL2}
  {\mathbb P} \bigl(|Z_{m}+ma| \le m\varepsilon_{m}+h(x) \ \ \forall  m
  \bigr)\to 1, \quad x\to\infty.
  \end{equation}
Based on Theorem \ref{Tnew2} and on \eqref{SLL2}, we obtain the following auxiliary
result
(see, e.g., Corollary 5 in \cite{SFSZ} 
for an analogous statement in the case of an ordinary random walk).

\begin{proposition}\label{prop1}
  Assume that the conditions of Theorem \ref{Tnew2} hold. Assume that $C>0$. 
  Let the function $h(x)\uparrow\infty$, $h(x)=\oh(x)$ be such that $F^I$ is
  $h$-insensitive, and introduce the events:
\begin{eqnarray*}
K_{n,x} &=& \bigcap_{m\le n-1}
\bigl\{ |Z_{m}+ma| \le m\varepsilon_{m}+h(x)\bigr\};\\
A_{n,x} &=& \{\xi_{n}>x+na\};\ \ 
A_{n,x}^{\varepsilon,h} = \bigl\{\xi_{n}>x+na+n\varepsilon_n+h(x)\bigr\}.
\end{eqnarray*}
  Then the following equivalences hold:
  \begin{eqnarray}
    \{M^{rw}>x\} &\sim & \label{3.8}
    \bigcup_{n\ge 1}\{M^{rw}>x\}\cap A_{n,x} \cap K_{n,x}
\ \sim\  \bigcup_{n\ge 1}\{M^{rw}>x\}\cap A_{n,x}^{\varepsilon,h} \cap K_{n,x} \nonumber \\
    &\sim &\label{3.9}
    \bigcup_{n\ge 1}\{M^{rw}>x\}\cap A_{n,x}
\ \sim\ \bigcup_{n\ge 1}\{M^{rw}>x\}\cap A_{n,x}^{\varepsilon,h}
\nonumber \\
    &\sim &\label{3.10}
    \bigcup_{n\ge 1}A_{n,x} \sim \bigcup_{n\ge 1} A_{n,x}^{\varepsilon,h}
    \end{eqnarray}
    and, therefore, 
    \begin{equation}\label{eq1000}
    {\mathbb P} (M^{rw}>x) 
\sim \sum_{n\ge 1}
{\mathbb P} \left( A_{n,x}^{\varepsilon,h}\cap K_{n,x}\right) 
\sim 
\sum_{n\ge
1}{\mathbb P} (A_{n,x}^{\varepsilon,h})
\sim 
\sum_{n\ge
1}{\mathbb P} (A_{n,x})
    \sim \frac{C}{a}\overline{F^{I}}(x).
    \end{equation}
    Finally 
    there exists a function $N=N(x)\to\infty$
    such that $\overline{F}^I(x+aN) \sim \overline{F}^I(x)$ and
    equivalences 
    \eqref{3.10} continue to hold if one replaces $n\ge 1$ by
    $n\ge N$.  
  \end{proposition}
\noindent {\em Proof.} One can easily verify
that 
$$
\bigcup_{n\ge 1} K_{n,x}\cap A_{n,x}^{\varepsilon,h}
\subseteq \{M^{rw}>x\}.
$$
The events $K_{n,x}\cap A_{n,x}^{\varepsilon,h}$ are disjoint and
$\sum_{n\ge 1} {\mathbb P} (A_{n,x}^{\varepsilon,h}\setminus A_{n,x}) =
o (\overline{F^I}(x))$, so
 by \eqref{SLL2},
 $$
 {\mathbb P} \bigl(
\bigcup_{n\ge 1} K_{n,x}\cap A_{n,x}^{\varepsilon,h}
 \bigr) =
 \sum_{n\ge 1} {\mathbb P}
 (K_{n,x}\cap A_{n,x}^{\varepsilon,h}) \sim
\sum_{n\ge 1} {\mathbb P}
 (K_{n,x}\cap {A}_{n,x}) \sim
\sum_{n\ge 1} {\mathbb P}
 ( {A}_{n,x}).
 $$
 Since ${\mathbb P} (M^{rw}>x) \sim \frac{C}{a}\overline{F^I}(x)$ by Theorem
 \ref{Tnew2} and
 since, by direct computations, $\sum_{n\ge 1} {\mathbb P}
 ( {A}_{n,x}) \sim \frac{C}{a}\overline{F^I}(x)$,
equivalences \eqref{eq1000} follow. 
The last fact follows directly from  Remark \ref{rem11}
and equivalences \eqref{3.10} and \eqref{eq1000}. \halmos

\vspace{0.3cm}

A special case of a modulated random walk is an ordinary
random walk with i.i.d. increments. 
Consider an auxiliary i.i.d. sequence
$\{{\xi}^{\sharp}_n\}$ with distribution $F$ and introduce the events
$$
A_{n,x}^{\sharp} = \{ {\xi}^{\sharp}_n>x+na \} 
\quad \mbox{and} \quad
D^{\sharp}_{x}=  \bigcup_{n\ge 1} A_{n,x}^{\sharp}.
$$
Assume there exists a function $e(x)\uparrow\infty$ such that, for
any $t>0$, there exists a limit
\begin{equation}\label{conv*}
\lim_{x\to\infty} \frac{{\mathbb
P}(D^{\sharp}_{x+te(x)})}{{\mathbb P}(D^{\sharp}_{x})}\, =\,
\overline{G}(t)
\end{equation}
with $\lim_{t\to\infty} \overline{G}(t)=0$. 
Remark that condition \eqref{conv*} is nothing else than condition
\eqref{conv} since ${\mathbb P} (D^{\sharp}_x) \sim \sum_{n\ge 1} {\mathbb P} (A_{n,x}^{\sharp})
\sim \frac{1}{a} \overline{F}^I(x)$.

On the event
$D^{\sharp}_x$, introduce the random variable\ 
$$
\tau^{\sharp} \equiv \tau^{\sharp}(x) = \min \{ n\ge 1 \ : \ {\mathbf 1}
(A_{n,x}^{\sharp}=1)\}.
$$
Then the following result holds:  
\begin{lemma}\label{cor3} 
Assume that the distribution $F^I$ is subexponential and that \eqref{conv}
holds. Then  
the conditional distribution of $\tau^{\sharp}/e(x)$
given 
$D^{\sharp}_x$ converges to the distribution $G$
(say, of the random variable $W$).
\end{lemma}
\noindent Indeed,
\begin{eqnarray*}\lefteqn{
\Prob(a\tau^{\sharp}/e(x) >t \ | \ \tau^{\sharp}<\infty ) \ =\ 
\Prob(\tau^{\sharp}>\frac{t}{a}e(x)  \ | \ \tau^{\sharp}<\infty )} \\
&\sim & \frac{\sum_{n> \frac{t}{a}e(x)}\Prob(\xi_n^{\sharp}>x+na)}{\Prob(
D^{\sharp}_x)} \ \sim\ 
\frac{\Prob(D^{\sharp}_{x+te(x)})}{\Prob(D^{\sharp}_{x})} 
\to  \overline{G}(t).
\end{eqnarray*}
\halmos

We now return to the modulated random walk.
On the event $\{ M^{rw}>x\}$, we similarly introduce the random variable
$$
\tau^{rw} = \tau^{rw} (x) = \min \{ n\ge 1 \ : \ Z_n>x \}.
$$
Recall from Proposition \ref{prop1}
that $$
\{ M^{rw}>x\}\sim D_{x}=  \bigcup_{n\ge 1} A_{n,x}.
$$
Then, by Lemma \ref{cor3}, 
we obtain: 
\begin{lemma}\label{cor4}
Under the assumptions of Theorem {\rm \ref{Tnew2}} with $C>0$ and \eqref{conv},
the conditional distribution of ${a\tau^{rw}}/e(x)$, conditioned on
$\{M^{rw}>x\}$, converges to the distribution $G$.
\end{lemma}
\noindent Indeed, the equivalence
$$
{\mathbb P} (a\tau^{\sharp} >te(x) \ | \
\tau^{\sharp}<\infty ) \sim {\mathbb P} (a\tau^{rw} >te(x) \
| \ \tau^{rw}<\infty )
$$
holds since we may represent conditional probabilities as ratios
of probabilities where both numerators and both denominators are
pairwise asymptotically proportional, with the same coefficient
$C$. 

Further, by Proposition \ref{prop1} and Lemma \ref{cor4},
one may deduce the following result.

\begin{theorem}\label{AK_ext}
Assume \eqref{conv} to hold. Then, under the conditions of Theorem
{\rm \ref{Tnew2}} and the assumption $C>0$, the distribution of
\begin{equation}\label{four}
\left(\frac{a\tau^{rw} }{e(x)}, \frac{Z_{\tau^{rw} -1}}{e(x)},
\max_{0\le m \le \tau^{rw} -1}
\frac{|Z_m+ma|}{\tau^{rw} }, \frac{Z_{\tau^{rw} }-x}{e(x)} \right),
\end{equation}
 conditioned on $\{M^{rw}>x\}$,
converges to the distribution of $(W, -W,0,W^{'})$
where $W$ and $W^{'}$ have the same distribution $G$ 
and, for
any positive $u$ and $v$,
\begin{equation}\label{psipsi}
{\mathbb P} (W >u, W^{'}>v ) = {\mathbb P} (W > u+v).
\end{equation}
\end{theorem}
\noindent This result is a complete analogue of Theorem 1.1 from
 Asmussen \& Kl\"uppelberg (1996) which was obtained in the case of
 an ordinary random walk.\\[1mm]
{\em Proof.} We have already proved the convergence of the first component in
\eqref{four}. From that and from \eqref{SLL2}, one may conclude
that
$$
\Prob \bigl(|Z_m+ma|\le m\varepsilon_m + h(x) \ \forall 
m<\tau^{rw}  \ | \ \tau^{rw}  < \infty \bigr) \to 1, \quad \mbox{as} \quad
x\to\infty.
$$
Then the convergence of the second and third components in
\eqref{four} follows if we take $h(x)\to\infty$ such that
$h(x)=\oh\bigl(e(x)\bigr)$.

It remains to
show the convergence of the last component in \eqref{four}.
Since  
$$
\bigl\{Z_{\tau^{rw} }-x > ve(x)\bigr\} \sim
\bigcup_{n\ge 1} \bigl\{ \xi_n-x-na > ve(x)\bigr\}, 
$$
we get
\begin{eqnarray*}
&& \Prob\bigl(Z_{\tau^{w} }-x > ve(x)\bigr) \sim
\Prob\biggl(\bigcup_{n\ge 1} \bigl\{ \xi_n-x-na > ve(x)\bigr\}\biggr)\\
&\sim & \sum_{n\ge 1} \Prob\bigl(\xi_n-x-na > ve(x)\bigr) =
\sum_{n\ge 1} \Prob\bigl(\xi_n >x + ve(x) +na \bigr)
\sim \Prob({D}_{x+ve(x)})
\end{eqnarray*}
(here we assume that $Z_{\tau^{rw}} =-\infty$ if $\tau^{rw}  = \infty$).
Similarly, equality \eqref{psipsi} follows since
\begin{eqnarray*}
\bigl\{a\tau^{rw}  > ue(x), Z_{\tau^{rw} }-x > ve(x)\bigr\}
&\sim & \bigcup_{n>ue(x)} \bigl\{ \xi_n >x+na,
\xi_n-na -x > ve(x) \bigr\} \\
&=&
\bigcup_{n>ue(x)} \bigl\{\xi_n > x+na+ve(x)\bigr\},
\end{eqnarray*}
and then
\begin{eqnarray*}\lefteqn{
\Prob\bigl(a\tau^{rw}  > ue(x), Z_{\tau^{rw} }-x > ve(x)\bigr)
\sim  \sum_{n>ue(x)}\Prob\bigl(\xi_n> x+na+ve(x)\bigr)}\\
&\sim & C\!\!\sum_{n>ue(x)}\Prob\bigl({\xi}^{\sharp}_n
>x+na+ve(x)\bigr)
\ = \ C\sum_{n\ge 1}\Prob\bigl({\xi}^{\sharp}_n
>x+na+(v+u)e(x)\bigr)\\
&\sim & C\, \Prob(D^{\sharp}_{x+(u+v)e(x)})
\ \sim\   \Prob(D_{x+(u+v)e(x)}).
\end{eqnarray*}
\halmos

\section{Continuous-time modulated regenerative processes}\label{S:ModRP}

We consider now a continuous-time process $Z(t)$ introduced in Section~\ref{S:Intr}
and assume that, more generally, it is a {\it regenerative process
which is modulated by a
 discrete-time regenerative process} $Y$. This means that
 (compare with the previous Section!)\\[1mm]
 (i) conditionally on $Y$, the random elements $V_{k+1}=\{Z(t)-Z(T_k), 0\le t \le R_{k+1}\}$
 are independent;\\
 (ii) for any $n$,
 \begin{equation}\label{VV}
 \Prob (V_n \in \cdot \ | \ Y) =\Prob (V_n \in \cdot \ | \ Y_n) \quad
 \mbox{a.s.}
 \end{equation}
 Let further, as in Theorem \ref{Th_ASS},
 $$
 \xi_k\ =\ Z(T_{k+1})-Z(T_k)\,,\ \ \
\xi^*_k\ =\ \sup_{T_k\le t<T_{k+1}}Z(t)-Z(T_k)
$$
and assume the conditions of Theorem \ref{Tnew2} and \eqref{conv} to hold. Then the statements of
Theorems \ref{Tnew2} and \ref{AK_ext} hold too. 

Note that $M \equiv\sup_{t\ge 0} Z(t)$ may be also represented as
$M= \sup_{n\ge 0} (\xi_1+\cdots+\xi_n+\xi^*_{n+1})$. Then we have the following result:
\begin{theorem}\label{MM}
Assume that the conditions of Theorem {\rm \ref{Tnew2}} and \eqref{conv} hold, and that $C>0$ in
Theorem {\rm \ref{Tnew2}}. Assume further that,
for all $y\in {\cal Y}$,
$$
{\mathbb P} (\xi^*_n >x \ | \ Y=y) \sim \overline{F}_y(x) \quad \mbox{a.s.}
$$
and that
$$
{\mathbb P} (\xi^*_n>x)\le c\overline{F}(x),
$$
for some $c\ge 1$ and all $x$.
Then, as $x\to\infty$,  
\begin{equation}\label{conti}
\{ M>x \} \sim \{ M^{rw}>x \}
\sim \bigcap_{n\ge 1} K_{n,x}\cap A_{n,x} 
\end{equation}
and, for ${\widehat{\tau}}^{rw} \equiv {\widehat{\tau}}^{rw}(x) = \min \{ n\ge 1 \ : \ Z_{n-1}+\xi^{*}_n>x\}$,
\begin{equation}\label{conti2}
{\mathbb P} (\tau^{rw}  = {\widehat{\tau}}^{rw} \ | \ M>x) \to 1, \quad x\to\infty 
\end{equation}
and
\begin{equation}\label{conti3}
{\mathbb P} (\tau^{rw}  = {\widehat{\tau}}^{rw} \ | \ M^{rw}>x) \to 1, \quad x\to\infty . 
\end{equation}

Therefore the statement of 
Theorem \ref{AK_ext} 
continues to hold if one
replaces in \eqref{four} $\tau^{rw} $ by ${\widehat{\tau}}^{rw}$ and then 
$Z_{\tau^{rw} }$ by $Z_{{\widehat{\tau}}^{rw}-1}+\xi_{{\widehat{\tau}}^{rw}}^{*}$.
\end{theorem}

The proof of Theorem \ref{MM} follows from routine minor modification of calculations
from the previous Section.

\section{Proof of Theorem~\ref{AFMain}}\label{S:RWProof}

Now we assume that the process $Z(t)$ is regenerative. This means that
$Y$ is a constant and, as a corollary, that conditions (3.1)-(3.4)
and (C1)-(C4) are redundant. Also, the $\xi_n$ are i.i.d. in this case and,
therefore, we may take $\tau^{\sharp}= \min \{n \ : \ \xi_n > x+na\}$. 

Let $\Probx$ denote the conditional probability given $\rtreg<\infty$,
write $T_n=R_1+\cdots+R_n$ and recall the definition of ${\widehat{\tau}}^{rw}$
from Theorem~\ref{MM}. Note that since
the events $\rtreg <\infty$ and 
${\widehat{\tau}}^{rw}<\infty$ 
coincide, and are
equivalent to each of the events $\rtrw<\infty$ and $\tau^{\sharp}<\infty$
(see Lemma \ref{cor4}), 
 we may use either of the four in conditioning arguments.

The proof of Theorem~\ref{AFMain} is a straightforward combination
of Theorem~\ref{AK_ext}, Theorem~\ref{MM} and
of the following two 
lemmas. 
Both use the fact, implicit in \cite{ASS} and also a consequence of \eqref{conti2}
and \eqref{conti3} of Theorem~\ref{MM}, 
that 
\begin{equation}\label{S1706a}
\Probx\bigl(\rtreg\in [T_{\rtrw-1},T_{\rtrw})\bigr)\ \to\ 1\,,\ \ x\to\infty
\end{equation}
since asymptotically
$$
\{\tau   \in [T_{\rtrw-1},T_{\rtrw})\}
\subseteq
\{ \tau^{rw} = {\widehat{\tau}}^{rw}\}.
$$

\begin{lemma}\label{lemmaSA0904}
$T_{\rtrw-1}/e(x)\,\to\, \mu W/a$ 
in $\Probx$-distribution. 
\end{lemma}
\noindent {\em Proof.} We use the representation
\begin{equation}\label{fra}
\frac{T_{\rtrw-1}}{e(x)} =
\frac{T_{\rtrw-1}}{\rtrw} \cdot \frac{\rtrw}{e(x)}.
\end{equation}
Choose $N=N(x)\to\infty$ from Proposition~\ref{prop1}. 
The first fraction in the LHS of \eqref{fra} converges to $\mu$ in ${\mathbb P}^{(x)}$ probability since,  
by the independence of $A_{n,x}$ and $T_{n-1}$ and by the SLLN, 
\begin{eqnarray*}\lefteqn{
\bigl\{\bigl|T_{\rtrw-1}/\rtrw-\mu\bigr|\le\varepsilon,\,\rtrw <\infty\bigr\}\ \sim\
\bigcup_{n\ge 1} \bigl\{\bigl|T_{n-1}/n-\mu\bigr|\le\varepsilon\bigr\}\cap A_{n,x}}\\ &\sim&
\bigcup_{n\ge N} \bigl\{\bigl|T_{n-1}/n-\mu\bigr|\le\varepsilon\bigr\}\cap A_{n,x}\ \sim\ 
\bigcup_{n\ge N} A_{n,x} \sim\ 
\bigcup_{n\ge 1} A_{n,x}\ \sim\ \{\rtrw<\infty\}.
\end{eqnarray*}
Then the second fraction
converges to $W/a$ by Theorem~\ref{AK_ext}, and the result follows.\halmos

\smallskip 

Recall that $\tau  = \sum_1^{{\widehat{\tau}}^{rw}-1} R_i + t_{{\widehat{\tau}}^{rw}}$.

\begin{lemma}\label{lemmaSA1706}
Under the conditions of Theorem~{\rm \ref{AFMain}}, 
$t_{{\widehat{\tau}}^{rw}}/e(x)\,\to 0\,$  
in $\Probx$-probability.
\end{lemma}
\noindent {\em Proof.}  %
%
By Theorem~\ref{AK_ext} and Theorem~\ref{MM}, for any $\delta >0$, one can choose $K>0$ such that
$$
\Prob\bigl({\widehat{\tau}}^{rw}/e(x) > K\bigr)\ \le\ \delta/2,
$$
for all $x$ large enough. Then, for any $y>0$, 
$$
\Probx\bigl(t_{{\widehat{\tau}}^{rw}}/e(x)>y\bigr)\ \le\ 
\bigl(1+\oh(1)\bigr) \frac{\sum_{n\le Ke(x)} \Prob\bigl(t_n>ye(x), \xi_n >x+na\bigr)}{
\overline{F}^I(x)/a} +\delta 
$$
If \eqref{weakSN} holds, then also $\limsup_{x\to\infty} e\bigl(x+ke(x)\bigr)/e(x) <\infty$,
for any $k>0$. Therefore,
the latter sum in the numerator is equivalent to
$$
\sum_{n\le Ke(x)}\Prob\bigl(t_n>y e(x+na) \, \big|\,  \xi_n > x+na\bigr) \Prob(\xi_n >x+na)
\ =\  \oh(1) \overline{F}^I(x)
$$
since the $(t_n,\xi_n)$ are i.i.d. 
Then we complete the proof by letting first $x\to \infty$ and then $\delta
\to 0$. \halmos

\smallskip

Combining this with the statements of Theorem~\ref{MM} and Lemma~\ref{lemmaSA0904}
completes the proof of Theorem \ref{AFMain}.

\section{Examples}\label{S:Ex}

\newtheorem{ex}[proposition]{Example}
\nc{\regpath}{{\cal P}}

\begin{ex}\label{ExSA21.6a}\rm
In the setting of Section \ref{S:RW}, we may assume that $Z(t)$ is a right-continuous
piecewise
constant process with $Z(n) = Z_n$. We shall show here that, under a natural extra assumption, 
Theorem~\ref{AFMain} holds for this model as well. Note that
because of the result of \cite{SFTKSZ}, we need not verify the
conditions of Theorem~\ref{Th_ASS} (which may be messy); all that
is needed is to establish \eqref{S2506c}.

Assume in addition that the distribution of the cycle length, $R$, has a lighter
tail than $\overline{F}(x)$, in the following strong sense: there exists  constant $c>1$
such that 
\begin{equation}\label{extra1}
\Prob(cR>x) = \oh\bigl(\overline{F}(x)\bigr), \quad x\to\infty.
\end{equation}
Let $\xi$ be the increment over the cycle, $\xi =\sum_1^R \xi_i$. 
For \eqref{S2506c} to hold, it suffices to show that, for any $y>0$, 
\begin{equation}\label{S2106a}
\Prob(R>yx,\,\xi>x)\ =\ \oh\bigl(\Fb(x)\bigr)\,,\ \ \ x\to\infty,
\end{equation}
where $F$ is the reference distribution. 
For any fixed $x_0$ and for $x\ge x_0$, as $x\to\infty$, 
\begin{eqnarray*}
\Prob(R>yx,\,\xi>x)
&=& \Exp (\Prob (\xi >x \ | \ R, Y_0,\ldots,Y_R) {\mathbf 1} (Y>yx))\\
&\le &
\sum_{k\ge xy} \overline{F^k}(x) \Prob (R=k) \quad \quad \mbox{(by property (C2))}\\
&\le &
\sum_{k\ge x_0y} \overline{F^k}(x) \Prob (R=k)\ \sim \
\overline{F}(x) \Exp [R; \,R>x_0y] 
\end{eqnarray*}
where, under assumption \eqref{extra1}, the last equivalence follows 
 from \cite{DDSFDK}, Theorem 1.
By letting $x_0\to \infty$, we obtain \eqref{S2106a}. \halmos 
\end{ex}

\begin{ex}\label{ExSA0201c}\rm
The Bj\"ork-Grandell model (\cite{a13}) is a regenerative risk process, such that
in addition to the cycle length $R$ also the rate $\Lambda$ of claims arrivals
within a cycle is random. All claims are i.i.d.\ with distribution $H$, with mean $m$, and independent
of $(R,\Lambda)$, and there is a constant rate $1$ of premium inflow. The infinite horizon ruin probabilities are discussed in \cite{a13}
for the light-tailed case and in \cite{ASS} for the heavy-tailed case. As noted in
\cite{ASS}, heavy tails of $\xi$ may occur in at least three ways: (i) $F$ is heavy-tailed;
(ii) $\Lambda$ is heavy-tailed; (iii) $R$ is heavy-tailed for sufficiently large values of
$\Lambda$.
Under some (not necessarily minimal) assumptions, we shall give arguments to identify
the limiting conditional behavior of $\rtreg$. 

For the following estimates, one may keep in mind that
\begin{equation}\label{S0401a}\xi\ =\ X_1+\cdots+X_{M_{R\Lambda}}-R\end{equation}
with $M$ an independent Poisson process at unit rate.
For the tail asymptotics of $\xi$, the $-R$ term may often be neglected
(see \cite{SARB} for some preliminary discussion and \cite{AAK} for a more complete picture).
Also, with light-tailed claims one may frequently approximate $X_1+\cdots+X_{M_{R\Lambda}}$
by $mR\Lambda$; the relevant large deviations arguments are given in detail in \cite{ASS} 
and will not be repeated here.

Consider first case (i) with  $R,\Lambda$ both light-tailed. Using \eqref{S0401a}
and an independence result from \cite{SARB}, it is standard that 
$$\Prob(\xi>x)\ \sim\ \Exp(R\Lambda)\Fb(x)\,.$$
By a classical inequality due to Kesten (see e.g. \cite{bingolteu87}, p.429), to each $\delta>0$ there is a $C_\delta<\infty$ such that
$\Prob(X_1+\cdots+X_n>x)\le C_\delta\e^{n\delta}\Fb(x)$ for all $n$. With 
$$p=\Prob(M_{R\Lambda}=1)>0\,,\quad q=p\Prob(X_1>x+R)\sim p\Fb(x)\,,$$ 
we get
\begin{eqnarray*}
\Exp[\e^{sR}\,|\,\xi>x]&=&\frac{\Exp[\e^{sR};\,\xi>x]}{\Prob(\xi>x)}\
\le\ \frac{1}{q}
\Exp\bigl[\e^{sR};\, X_1+\cdots+X_{M_{R\Lambda}}>x\bigr]\\
&\le&\frac{1}{q}
\Exp\bigl[\e^{sR}C_\delta\e^{\delta M_{R\Lambda}}\Fb(x)\bigr]
\ \sim\ \frac{C_\delta}{p}
\Exp\bigl[\e^{sR}\e^{R\Lambda(\e^\delta-1)}\bigr]\,.
\end{eqnarray*}
Taking $s,\delta$ small enough, this expression is finite, and its independence of $x$
together with $e(x)\to\infty$ then gives \eqref{S2506c} and the conclusion of
Theorem~\ref{AFMain}.

Consider next case (ii) with  $F$ light-tailed and $(R,\Lambda)$
satisfying $\Prob(\Lambda>x)\sim x^{-\alpha}$ with $\alpha>1$ and $\Exp R^{\alpha'}<\infty$ for some $\alpha'>\alpha$. Then by Breiman's theorem (\cite{Breiman}, \cite{ClineSam}, \cite{MaulZw}),
$\Prob(mR\Lambda>x)\sim cx^{-\alpha}$ where $c=m^\alpha\Exp R^\alpha$. 
By a large deviations argument,
$$\Prob(M_{R\Lambda}>x)\ =\ \Prob(R\Lambda>x)\,+\,\Oh(\e^{-\varepsilon_1 x})$$
for some $\varepsilon_1>0$.
A further large deviations argument given in \cite{AKS} then shows that
$$\bigl\{X_1+\cdots+X_{M_{R\Lambda}}>x\bigr\}\Delta \{mR\Lambda>x\}\ =\ A(x)$$
where $\Prob A(x)=\Oh(\e^{-\varepsilon x})$ for some $\varepsilon>0$ .
In particular $X_1+\cdots+X_{M_{R\Lambda}}$ 
has asymptotic tail $cx^{-\alpha}$. Hence so has $\xi=$ 
$X_1+\cdots+X_{M_{R\Lambda}}-R$ (see \cite{SARB}; note that this is non-trivial due to dependence).
Let $\alpha<\alpha''<\alpha'$ and let $R^*$ 
be a r.v.\ with distribution
$$\Prob(R^*\in \dd t)\ =\ \Exp\bigl[R^{\alpha''};\,
R \in\dd t\bigr]/\Exp(R^{\alpha''})\,.$$
Then
\begin{eqnarray*}
\Exp[R^{\alpha''} \ | \ \xi >x ]&\sim&
\frac{\Exp[R^{\alpha''};\,\xi>x]}{cx^{-\alpha}}\ 
\le\ \frac{1}{cx^{-\alpha}}
\Exp\bigl[R^{\alpha''};\,X_1+\cdots+X_{M_{R\Lambda}}>x\bigr]\\
&=&\frac{1}{cx^{-\alpha}}
\Exp\bigl[R^{\alpha''};\,mR\Lambda>x\bigr]\,+\,\Oh(\e^{-\varepsilon x})\\ &=&
\frac{1}{cx^{-\alpha}}\Exp R^{\alpha''}\Prob(R^*\Lambda>x/m)\,+\,\Oh(\e^{-\varepsilon x}).
\end{eqnarray*}
 Another application of Breiman's theorem justified by the choice of $\alpha''$ shows that
 $R^*\Lambda$ has a distribution tail asymptotically proportional to $x^{-\alpha}$. Hence
 $\Exp[R^{\alpha''}\,|\,\xi>x]$ stays bounded as $x\to\infty$, and arguing as above
 gives \eqref{S2506c} and the conclusion of Theorem~\ref{AFMain}.

In contrast, the behavior in case (iii) is different, see Section~\ref{S:NSGrowth}.
\halmos\end{ex}

\begin{ex}\label{ExSA0201a}\rm
Let $Z$ be a two-stage fluid model, where a cycle $R$ is composed of two stages such
that the first has deterministic length $a_1$ and the second a random length $R_2$ with
a subexponential distribution $F$ with mean $a_2<a_1$. In stage 1, $Z$ decreases deterministically at rate
1 and in stage 2, $Z$ increases deterministically at rate 1 (thus $-a=a_2-a_1<0$). Clearly,
$\xi>x$ occurs if and only $R_2>x+a_1$. Thus
$\Prob(\xi>x)=\Fb(x+a_1)\sim\Fb(x)$ and given $\xi>x$, $R$ is at least $x$. Since $e(x)=\Oh(x)$
in all examples, condition \eqref{S2506c} can not hold and more precisely, given $\xi>x$, $R$
is of order $x+e(x)$. Therefore $\rtreg$ is of order $e(x)$ in the regularly varying case
but with a larger multiplier than $\mu W$, and of order $x>\!\!>e(x)$ for other subexponential
distributions.

Note that this example shows that the regenerative setting is more flexible than
the Markov additive one: if one considers the discrete time analogue, the increments in each
Markov stage are bounded and there is thus no version of condition (C2) of Section~\ref{S:RW}
with $F$ heavy-tailed. On the other hand, conditions may be easier to verify in the
Markov additive setting.
\end{ex}\halmos

\section{Different  growth rates}\label{S:NSGrowth}

The following  result is straightforward given Lemma~\ref{lemmaSA0904}
and the proof of Theorem~\ref{AFMain}. For simplicity (to avoid distinction between
$x$ and $e(x)$) we state it only for the regularly varying case where the r.v.\ $W$ 
in Theorem~\ref{cor1} is Pareto:

\begin{cor}\label{Cor0904a} Assume that $F$ in \eqref{S15.6b} is regularly varying and that
instead of Condition~\eqref{S2506c} we have 
\begin{equation}\label{S0904c}
\Prob\bigl(t_1/e^*(x)>y\,\big|\,\xi>x) \ \to\ \Prob(W^*>y)\ \ \text{for all  }y,
\end{equation}
some function $e^*(x)$ with $\liminf e^*(x)/x>0$ and some r.v.\ $W^*$.\\
{\rm (i)} If $e^*(x)\sim dx$ for some $d$, then
$$\frac{\rtreg}{x}\ \to \ W\mu/a+d(1+W)W^*\ \ \text{in }\Probx-\text{distribution};$$
{\rm (ii)} if $e^*(x)/x\to\infty$, then
$$\frac{\rtreg}{e^*\bigl(x(1+W)\bigr)}\ \to \ W^*\ \ \text{in }\Probx-\text{distribution}$$
with $W,W^*$ independent in both {\rm (i)} and {\rm (ii)}, and $W$ independent of $\rtreg$
in {\rm (ii)}. In particular, if $e^*(x)\sim dx^\beta$ with $\beta> 1$, then
$\rtreg/x^\beta\,\to\,d(1+W)^\beta W^*$ in $\Probx$-distribution.
\end{cor}
\proof The asymptotic $\Probx$-distribution of $\tau$ is the same as the asymptotic
distribution of $\sum_1^{\widehat \tau^{\rm rw}-1}R_i+t_{\widehat \tau^{\rm rw}}$.
Here $\sum_1^{\widehat \tau^{\rm rw}-1}R_i/x\to \mu W/a$ 
in $\Probx$-distribution (Lemma~\ref{lemmaSA0904}).
More generally,
$$\frac{1}{x}\Bigl(\sum_1^{\widehat \tau^{\rm rw}-1}R_i,\,\sum_1^{\widehat \tau^{\rm rw}-1}\xi_i\Bigr)\ \to\ 
(\mu W/a,W)\,.$$
Given $W=w$, $\xi_{\widehat \tau^{\rm rw}}$ will asymptotically have to exceed
$xw+x$, implying $t_{\widehat t^{\rm rw}}/e^*\bigl(x(w+1)\bigr)\to W^*$ and the conclusion
of (i) since the limit $W^*$ does not depend of $w$. 
For (ii), just note that in this case $\sum_1^{\widehat \tau^{\rm rw}-1}R_i$ may be
neglected.\halmos

\smallskip

We next first give an example of $e^*(x)\sim dx$ and thereafter some discussion of 
what may happen if  $e^*(x)/x\to\infty$.

\begin{ex}\label{ExSA0201d}\rm
We return to the Bj\"ork-Grandell model 
in case (iii). Here one expects  that given
$\Lambda=\lambda$, the surplus process
$\sum_1^{N(t)}U_i-t$ can be approximated by $\lambda m t-t$, and this is confirmed by
the large deviations bounds in \cite{ASS}. Therefore the behavior should be
like a fluid model with heavy-tailed on periods, so that the exceedance time of $x$
within a cycle should be of order $x$ and accordingly makes a genuine contribution
to $\rtreg$.

We  next verify this  statement and and make it more precise, assuming as in \cite{ASS} that claims are light-tailed and independent
of $(R,\Lambda)$,
that for some $\lambda_0>1/\mu$
\begin{eqnarray*}
\Prob\bigl(R>t\,\big|\,\Lambda=\lambda\bigr)&=&\Fb(t)\,,\ \ \lambda>\lambda_0\,,\\
\Prob\bigl(R>t\,\big|\,\Lambda=\lambda\bigr)&\le&\Gb(t)\,,\ \ \lambda\le\lambda_0,
\end{eqnarray*}
for some regularly varying $F$ with $\Fb(t)={L(t)}/{t^\alpha}$ ($L$ slowly varying), some $G$
satisfying $\Gb(t)=\oh\bigl(\Fb(t)\bigr)$, and  the following regularity condition:
\begin{equation}\label{S0904f}
\sup_{x\ge x_0}\frac{L(x/y)}{L(x)}\ \le g(y)
\end{equation}
for all $y>0$, some $x_0>0$ and some function $g(y)$ with 
$\Exp\bigl[\Lambda^\alpha g(\Lambda)\bigr]<\infty$.

It is then shown in \cite{ASS} that the conditions of Theorem~\ref{Th_ASS} are satisfied and that
\begin{equation}\label{28.3a}
\psi(x)\,\sim\,c_1\Fb(x)\ \ \text{where}\ c_1\,=\,
\frac{c}%
{(\alpha-1)\bigl[\Exp R-m\Exp(\Lambda R)\bigr]}\,,\ \ c\,=\,\Exp\bigl[(\Lambda m-1)^\alpha;\,\Lambda>\lambda_0\bigr]\,.
\end{equation}
This depends on the estimate
\begin{eqnarray}\label{8.4a}
\Prob(\xi>x)&\sim&c\Fb(x)\,.
\end{eqnarray}

As preparation for the study of the ruin time, we first recall the proof of \eqref{8.4a}. That the event
$\xi>x$ occurs is by the LD arguments equivalent to $R(\Lambda m-1)>x$, and so
\begin{eqnarray}\label{8.4b}
\Prob(\xi>x)&\sim&\int_{\lambda_0}^\infty 
f_\Lambda(\lambda)\Fb\bigl(x/(\lambda m-1)\bigr)\,\dd\lambda \nonumber \\
&=&\int_{\lambda_0}^\infty 
f_\Lambda(\lambda)(\lambda m-1)^\alpha\frac{L\bigl(x/(\lambda m-1)\bigr)}{x^\alpha}\,\dd\lambda \nonumber\\
&\sim&\int_{\lambda_0}^\infty 
f_\Lambda(\lambda)(\lambda m-1)^\alpha\frac{L(x)}{x^\alpha}\,\dd\lambda\ = \ c\Fb(x)\,,
\end{eqnarray}
where the last $\sim$ follows by dominated convergence justified by \eqref{S0904f}.
If $\xi>x,\tau\le xt$ is to occur, we need in addition $x/(\lambda m-1)\le xt$, and so
by the same dominated convergence argument
\begin{eqnarray}\label{8.4c}
\Prob(\xi>x,\tau\le xt)&\sim&\int_{\lambda_0\vee (1/t+1)/m}^\infty 
f_\Lambda(\lambda)\Fb\bigl(x/(\lambda m-1)\bigr)\,\dd\lambda \nonumber \\
&\sim&\int_{\lambda_0\vee (1/t+1)/m}^\infty 
f_\Lambda(\lambda)(\lambda m-1)^\alpha\frac{L(x)}{x^\alpha}\,\dd\lambda\ = \ c\Fb(x)W^*(t)\,,
\end{eqnarray}
where $W^*$ is the distribution with c.d.f.\
$$\Prob(W^*\le t)\ =\ \Biggl\{\begin{array}{cl}\displaystyle \frac{1}{c}\int_{(1/t+1)/m}^\infty 
f_\Lambda(\lambda)(\Lambda m-1)^\alpha\,\dd\lambda , &t\le 1/(\lambda_0m-1)\\
1, &t> 1/(\lambda_0m-1)\end{array}\,.$$
From Corollary~\ref {Cor0904a} we therefore conclude that
$\rtreg(x)/x\,\to\,W^*(1+W)$ in $\Probx$-distribution.\halmos
\end{ex}

We proceed  to discussing when $\rtreg$ may grow
at larger rates than $e(x)$ and how fast the rate may be. 
If $t_x=\Exp[R\,|\,\xi=x]\to\infty$ faster than $e(x)$, one expects $R$ given $\xi>x$ (and hence
often $\rtreg$) to grow at a faster rate than $e(x)$.
At first sight, one could conjecture
that any rate is possible. This is, however, not possible because of the requirement $\Exp R
<\infty$. Suppose, for example, that $F$ is a discrete subexponential distribution
with point probabilities $f_x=\Prob(\xi=x)\sim c_1/x^{\alpha+1 }$.
Assuming $t_x\sim c_2 x^\beta$, we then get
$$\infty\ >\  \Exp R \ =\ \sum_0^\infty t_xf_x\ \approx \ \sum_0^\infty c_2 x^\beta\,c_1/x^{\alpha+1 }
\,,$$
implying $\beta<\alpha$. The following result gives the more precise upper bound
$c/\Fb(x)$ and is more satisfying by being in terms of the growth rate of $\rtreg$
rather than expected values:

\begin{theorem}\label{STh0201b} Let $F(x)=\Prob(\xi\le x)$ be a discrete subexponential distribution
with point probabilities $f_0,f_1,\ldots$ and $\varphi$ a function with 
$\varphi (x)/e(x)\to\infty$. 
Assume that $\Prob\bigl(R>\varepsilon \varphi (x)\,\big|\,\xi>x)$ $\ge \delta$ for some $\varepsilon,\delta>0$
and all large $x$.
Then $\varphi (x)\le c/\Fb(x)$ for some constant $c$.
\end{theorem}
\proof Define $t(x)$ as above and let 
\begin{equation}\label{S0201a}k(x)\ =\ \Exp[R\,|\,\xi>x]\ =\ 
\frac{1}{\Fb(x)}\bigl(t_{x+1}f_{x+1}+t_{x+2}f_{x+2}+\cdots\bigr)
\end{equation}
Multiplying by $\Fb(x)$ and subtracting the resulting equation with $x$ replaced by $x+1$, it follows that
$$
t_{x+1}\ =\ \frac{1}{f_{x+1}}\bigl(k(x)\Fb(x)-k(x+1)\Fb(x+1)\bigr)\,.
$$
This expression needs to be positive which gives $k(x+1)/k(x)<\Fb(x)/\Fb(x+1)$ and,
multiplying from $x=1$ to $y-1$,
\begin{equation}\label{S0201b}k(y)\ <\ k(1)\Fb(1)\frac{1}{\Fb(y)}\,.\end{equation}
However, clearly
$k(x)\,=\,\Exp[R\,|\,\xi>x]\ \ge\ \varepsilon\delta h(x)\,,
$
from which we conclude
$$\varphi (x)\ \le\ \frac{k(x)}{\varepsilon\delta}\ <\ \frac{k(1)\Fb(1)}{\varepsilon\delta}\,\frac{1}{\Fb(x)}
\eqno{\halmos}$$

\smallskip

That the upper bound  of order $1/\Fb(x)$ is attainable follows from the following example:
\begin{ex}\label{ExSA0201b}\rm
Let $F$ be a discrete subexponential distribution with point probabilities $f_0>0,f_1,f_2,\ldots$\ \,
A discrete-time regenerative process $Z$ is constructed as follows.
At the start of a cycle, a r.v.\ $X$ with distribution $F$ is drawn. If $X=0$, one takes 
$R=1,\xi=-b$. If $X=x>0$, one takes $R=\varphi (x)$ for some suitable 
$\varphi (x)\uparrow\infty$, 
$Z_0=\ldots Z_{\varphi (x)-2}=0$,
$\xi=Z_{\varphi (x)-1}=x$ (one then needs to choose $b$ such that $\Exp\xi<0$).

The question is whether all rates are attainable. To discuss this point, let 
$\varphi (x)\uparrow\infty$.
By Theorem~\ref{STh0201b}, $\varphi$ must satisfy $\varphi (x)
=\Oh\bigl(1/\Fb(x)\bigr)$. Conversely, the construction
works if \eqref{S0201b} holds and gives a risk process such that 
$\rtreg$ grows at rate at least $\varphi (x)$. 

For example, if $F$ is regularly varying with
index $\alpha>1$, this allows for growth rates $\varphi (x)$ 
of order $x^\beta$ with $1<\beta<\alpha$,
whereas $e(x)$ only is of order $x$.
\end{ex}


\end{document}